\documentclass[11pt]{amsart}
\newif\ifpdf
    \ifx\pdfoutput\undefined
    \pdffalse 
    \else
    \pdfoutput=1 
    \pdftrue
    \fi

    \ifpdf
    \usepackage[pdftex]{graphicx}
    \else
    \usepackage{graphicx}
    \fi
\usepackage{amsmath,amssymb,amsthm,epsfig}

\newcommand\Z{\mathbb Z}

\newtheorem{theorem}{Theorem}[section]
\newtheorem{proposition}[theorem]{Proposition}

\theoremstyle{definition}

\theoremstyle{example}

\input{epsf.tex}

\ifpdf
    \DeclareGraphicsExtensions{.pdf, .jpg, .tif}
    \else
    \DeclareGraphicsExtensions{.eps, .jpg}
    \fi

\title{Distortion of wreath products in some finitely presented groups}
\author{Sean Cleary}
\address{Department of Mathematics \\
The City College of New York  and the CUNY Graduate Center\\
City University of New York \\
New York, NY 10031}
\email{ cleary@sci.ccny.cuny.edu}
\urladdr{http://www.sci.ccny.cuny.edu/\~ cleary}
\thanks{The author gratefully acknowledges support from PSC-CUNY grant \#64459-0033
and the hospitality of the Centre de Recerca Matem\`atica }

\begin{document}

\begin{abstract}
Wreath products such as $\Z \wr \Z$ are not finitely-presentable yet
can occur as subgroups of finitely presented groups.  Here we compute
the distortion of $\Z \wr \Z$  as a subgroup of Thompson's group $F$ and 
as a subgroup of  Baumslag's metabelian group $G$. 
 We find that $\Z \wr \Z$  is undistorted in $F$ but
is at least exponentially distorted in $G$.
 \end{abstract}
 
\def\cprime{$'$}

\maketitle

\section{Introduction}
\label{sec:intro}
Here we consider aspects of the question of the distortion of infinitely-related
groups as subgroups of finitely-presented groups.
Higman \cite{higman} showed that every recursively-presentable group occurs as a subgroup of a finitely-presented group, but it is not clear in general what happens to the geometry of the group since
this embedding uses complicated algebraic methods  and methods from recursive function
theory which may affect the geometry of the group severely.  
Ol{\cprime}shanski{\u\i} \cite{subdistort} 
constructs isometric embeddings of recursively presentable groups into finitely presented groups
using difficult methods that do not lead to easily constructed examples.
 In the particular concrete cases here, we consider concrete embeddings of
 one of the simplest finitely-generated but not finitely-presentable
groups, $\Z \wr \Z$.  We consider two embeddings of $\Z \wr \Z$ into finitely presented groups.
The first is as a subgroup of Thompson's group $F$ and the second is as subgroup of
Baumslag's remarkable finitely presented metabelian group which contains $\Z \wr \Z$ and thus
a free abelian subgroup of infinite rank.  The distortion of the metric of $\Z \wr \Z$ is linear 
in Thompson's group $F$ but is exponential in Baumslag's group.

\section{Background}
\label{sec:background}

\subsection{Metrics of wreath products}

We construct the wreath product in the standard manner, as
special case of a semi-direct product. Given two groups $G$ and $H$,
we form the wreath product $G \wr H$ by taking the direct product
of $|H|$ copies of $G$ with copy of $G$ indexed by an element of $H$.
The generators of
$G$ act on the conjugate copies  of $G$, while
generators of $H$ act on the coordinates to determine to which of  these conjugate copies of
$G$ the generators of $G$ will be applied.

Two of the simplest infinite wreath products are the lamplighter
group $\Z_2 \wr \Z$ and
 $\Z \wr \Z$.  Cleary
 and Taback \cite{deadlamp} analyzed aspects of the metric geometry
 of those groups and other wreath products.
 There are  natural normal forms for elements
in these groups which lead to geodesic words for elements in these
groups with respect to their standard generating sets.

For  $\Z \wr \Z$, we consider the standard presentation:

$$<a,t | [a^{t^i},a^{t^j}],  \mbox{ for } i,j \in \Z>$$

where we denote the conjugate $b^{-1} a b$ of $a$ by $b$ as $a^b$ and
the commutator $a b a^{-1} b^{-1}$ by $[a,b]$.

\begin{figure}
\includegraphics[width=3.5in]{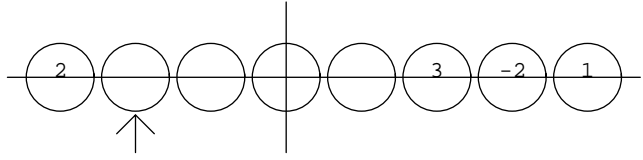}
\caption{Diagram for $w=a_2^3 a_3^{-2} a_4 a_{-3}^2 t^{-2}$.  The origin in the
wreath product direction is indicated by a vertical line, empty circles denote counters
which are zero, and the final cursor position is indicated by the arrow. \label{fig:wreathpic}}
\end{figure}

Geometrically, we can think about this wreath product as a 
set of parallel copies of $\Z$ strung together along their respective origins.
We can think of this as a string of counters, arranged from left to right with
one counter distinguished as the origin.    As in the lamplighter group, we
imagine a cursor which moves along the string of counters and will
point to a particular one of these counters
as being of current interest.  The generator $a$ acts as a generator of $\Z$
in the factor to which the cursor currently points and increases the counter
in that factor, and the generator $t$ moves
the cursor to the right to the next counter. A typical such word is illustrated in
Figure \ref{fig:wreathpic}. 

The starting configuration of these counters, corresponding to the identity
element in $\Z \wr \Z$, is with all of the counters at zero and the cursor 
resting at the counter designated at the origin.
We consider a word in these generators as a sequence
of instructions to move the cursor and change the counter in the current factor.
After application of a long string of the generators, we will  be
in a state where a finite number of counters are non-zero and the cursor
points at a particular counter, called the ``final position'' of the cursor for that word.

We define $a_n= a^{t^n}$ and note that $a_n$ is a generator of the conjugate
copy of $\Z$ indexed by $n$.  These $a_n$ commute and we can put any word
in the generators into one of two normal forms, `right-first' and `left-first', as described
by Cleary and Taback in \cite{deadlamp}:

$$rf(w)=a^{e_{1}}_{i_1} a^{e_{2}}_{i_2} \ldots a^{e_{k}}_{i_k} a^{f_{1}}_{-j_1} a^{f_{2}}_{-j_2} \ldots a^{f_{l}}_{-j_l} t^{m}$$
or
$$lf(w)=a^{f_{1}}_{-j_1} a^{f_{2}}_{-j_2} \ldots a^{f_{l}}_{-j_l} a^{e_{i}}_{i_1} a^{e_{i}}_{i_2} \ldots a^{e_{k}}_{i_k}   t^{m}
$$
with $ i_k > \ldots i_2 > i_1 \geq 0 $ and $ j_l > \ldots j_2 >
j_1 > 0 $ and $e_i, f_j \neq 0$.

The final resting position of the cursor is easily seen to be $m$ from either of
these normal forms, and we can see that the leftmost non-zero counter is in position
$-j_l$ and the rightmost non-zero counter is in position $i_k$.

In the `right-first' form, $rf(w)$, the cursor moves first to the
right from the origin, changing the counters in the appropriate
factors as the cursor moves to the right.   Then the cursor
moves back to the origin not affecting any of the counters until passing
the origin.  Past the origin, the cursor continues to works leftwards, again changing
the counters in the appropriate factors.  Finally, the cursor moves to
its ending location from the leftmost nonzero counter to the left of the origin.

The `left-first' form is similar, but instead of initially moving
to the right, the cursor begins by moving toward the left.

At least one of these normal forms will lead to
minimal-length representation for $w$, depending upon the final
location of the cursor.  If $m$ is non-negative, then the left-first
normal form will lead to a geodesic representative, and if $m$ is
non-positive, the right-first normal form will lead to a geodesic
representative, as described in \cite{deadlamp}, which gives
the following measurement of length:

\begin{proposition}[Prop 3.8 of \cite{deadlamp}]
If a word $w \in \Z \wr \Z$ is in either normal form given above, we measure
the word length of $w$ with respect to $\{a,t\}$ and have
$$|w|=\sum_{n=1}^{k} |e_{i_n}|+\sum_{n=1}^{l} |f_{j_n}| + min\{2 j_l+i_k + | m-i_k|, 2 i_k+j_l+|m+j_l|\}.$$
\end{proposition}

The first two terms are the minimum number of applications of $a^{\pm1}$ needed to put
all of the counters into their desired states and the last term 
is the  minimum possible amount of movement required to visit the
left- and right-most non-zero counters and then the final position of the cursor, counting
the required applications of $t^{\pm1}$.

The word $w=a_2^3 a_3^{-2} a_4 a_{-3}^2 t^{-2} $ pictured  in
Figure \ref{fig:wreathpic} has geodesic representatives in right-first normal form
since the final position of the cursor is to the left of the origin.  One such
minimal representative is $ t^2 a^3 t a^{-2} t a t^{-7} a^2 t$ of length 20.

\section{$\Z \wr \Z$ as a subgroup of Thompson's group $F$}
\label{sec:thomp}

Thompson's group $F$ is a remarkable, finitely-generated, finitely-presented
group which can be understood via a wide range of perspectives.
Cannon, Floyd, and Parry \cite{cfp} give an excellent overview of
the properties of $F$.   The standard infinite presentation of $F$ is
given by:

$$<x_0, x_1, \ldots | x_n^{x_i}= x_{n+1} \mbox{ for } i < n >.$$

Since $x_2=x_1^{x_0}$ and so on, $F$ is generated by the first two generators
and we can define $x_{n+1}= x_n^{x_0}$ to express all generators and thus
all group elements in terms of $x_0$ and $x_1$. 
Furthermore, all of these infinitely many relations are consequences of the first
two non-trivial relations, so we have the standard finite presentation:
$$<x_0, x_1 | x_2^{x_1}= x_3, x_3^{x_1}= x_4>.$$

Thompson's group $F$ can be described in terms of rooted tree pair diagrams,
and there is a natural method of converting between words in a normal form
with respect to the infinite generating set and tree pair diagrams, via the
method of leaf exponents, 
as described in Cannon, Floyd and Parry \cite{cfp}.  There is a natural notion
of a reduced tree pair diagram described there and there are efficient means
to convert between the unique normal form for an element of $F$ and the
unique reduced tree pair diagram for that word. 

We consider a rooted binary tree with $n$ leaves as being constructed of
$n-1$``carets,'' which are interior nodes of the tree together with the two downward
directed edges from that node.  The ``left side'' of a tree consists of nodes and
edges
which are connected to the root by a path consisting only of left edges, and
similarly the ``right side'' of a tree consists of nodes and
edges
which are connected to the root by a path consisting only of right edges.
A tree pair diagram $(S,T)$ is made up of a `positive' tree $T$ and a 
`negative' tree $S$.

\begin{figure}
\includegraphics[width=3.5in]{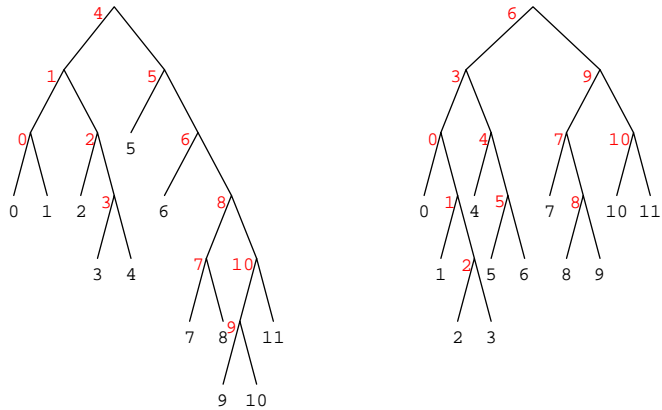}
\caption{Tree pair diagram for
 $x_0^2 x_1 x_2 x_4 x_5 x_7 x_8 x_9^{-1} x_7^{-1} x_3^{-1} x_2^{-1} x_0^{-2}$
with carets and leaves numbered. \label{fig:treepairexample}}
\end{figure}

\begin{figure}\includegraphics[width=3.5in]{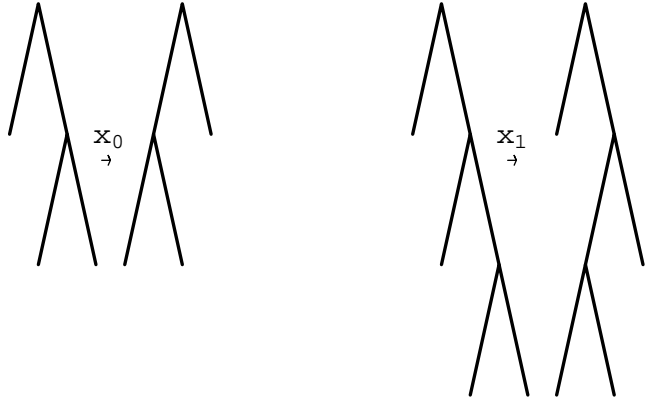}\\
\caption{The tree pair diagrams representing the generators $x_0$
and $x_1$ of $F$. \label{fig:x0x1}}
\end{figure}

The reduced tree pair diagrams for $x_0$ and $x_1$ and a typical word
in normal form are pictured in Figures  \ref{fig:treepairexample} and \ref{fig:x0x1}.

To understand the metric properties of $F$, we consider expressing words with respect to
the finite generating set.     Burillo, Cleary and Stein \cite{bcs} estimated the word length
in terms of the number of carets and showed that the number of carets is quasi-isometric
to the word length.  Fordham \cite{blakegd} developed a remarkable method using tree pair diagrams
to efficiently compute exact word length and find minimal length representatives of words.

We can understand word length of elements represented as tree pair diagrams  by
understanding the how the generators change the tree pair diagram for $w$ to that for $wg$ for the
generators, as described in Fordham \cite{blakegd} and in Cleary and Taback \cite{ct2}.
The right actions of the generators can be described as `rotations' which change the
tree.

\begin{figure}
\includegraphics[width=3.5in]{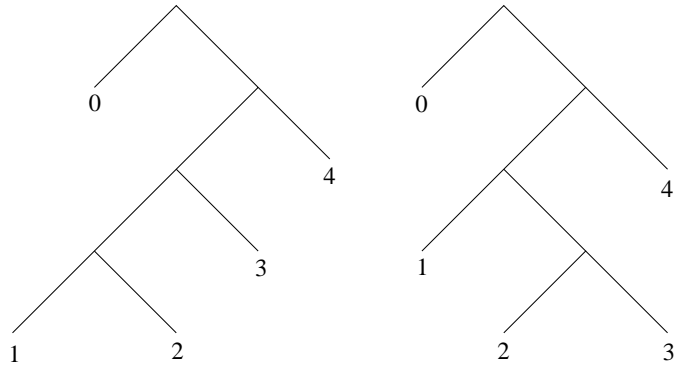}
\caption{Tree pair diagram for
 $ x_1  x_2 x_1^{-2}$, the image of $a$ under $\phi$. \label{fig:h}}
\end{figure}

The wreath product $\Z \wr \Z$ is a subgroup of $F$ and can be realized in many different
ways.  Perhaps the simplest is as the subgroup $H$ generated by $x_0$ and $h= x_1  x_2 x_1^{-2}$,
pictured in Figure \ref{fig:h}.  The isomorphism between this subgroup and $\Z \wr \Z$ is
given by the homomorphism $\phi(w):L \rightarrow F$ where $\phi(a)= x_1  x_2 x_1^{-2}$
and $\phi(t)=x_0$.   The isomorphism is readily established after we see that
$x_0$ conjugates $\phi(t)$ and its powers to elements that
commute and that there are no other relations.

To understand the distortion of the subgroup $H$ in $F$, we compare the
word length of an element $w= a^{e_{1}}_{i_1} a^{e_{2}}_{i_2} \ldots a^{e_{k}}_{i_k} a^{f_{1}}_{-j_1} a^{f_{2}}_{-j_2} \ldots a^{f_{l}}_{-j_l} t^{m}$ with its image in $F$. 

\begin{theorem}
The subgroup $H$ isomorphic to  $\Z \wr \Z$  in $F$ generated by $x_0=\phi(t)$ and $h= x_1^2 x_2^{-1} x_1^{-1}=\phi(a)$
is undistorted.
\end{theorem}

{\bf Proof:}
We count the number of carets of the image of a word 
$w$. 
First, we consider the case when $m=0$ and then the  cases where $m$ is nonzero.

\begin{figure}\includegraphics[width=5.5in]{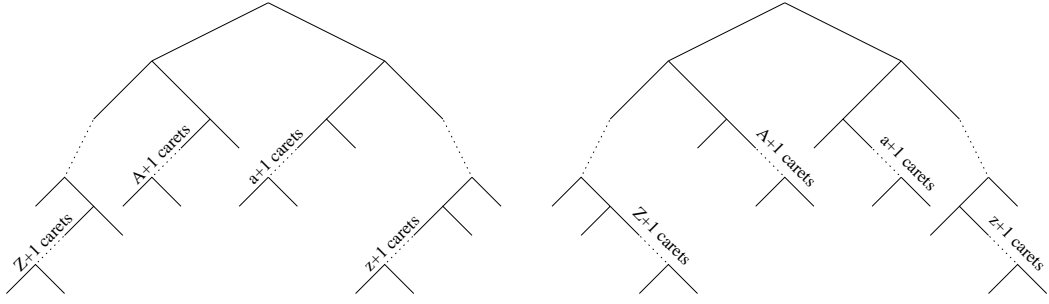}\\
\caption{The tree pair diagram for  the image $\phi(w)$ of a word $w=a_0^a a_1^b \ldots a_n^z a_{-1}^A a_{-2}^B \ldots a_{-m}^Z$ with $t$ exponent sum 0  and all positive exponents for $a_i$.\label{fig:typicalm}}
\end{figure}

{\bf Case $m=0$:}
Here, the image of the word as a tree pair diagram has a characteristic form where the
root of the positive tree is paired with the root of the negative tree, such as that shown in
 Figure \ref{fig:typicalm}.  In the general case where both $k$ and $l$ are positive, we
have the following carets:
\begin{itemize}
\item A single root caret
\item $i_k+1$ right carets
\item $\displaystyle \sum_{n=1}^k (| e_n|+1) $ interior carets below the right arm of the tree
\item $j_l$ left carets
\item $\displaystyle \sum_{n=1}^l (| f_n|+1) $ interior carets below the left arm of the tree
\end{itemize}

This gives a total  $N(\phi(w))=i_k+j_l+2+ \sum_{n=1}^k (| e_n|+1) + \sum_{n=1}^l (| f_n|+1) $ carets
in the image of $w$.
By Burillo, Cleary and Stein \cite{bcs}, the number of carets is quasi-isometric to the
word length in $F$ with respect to $\{x_0,x_1\}$ and since the length of $w$ in 
  $\Z \wr \Z$ is $ 2 j_l+2 i_k +\sum_{n=1}^{k} |e_{n}|+\sum_{n=1}^{l} |f_{n}| $,
  we see that these lengths are quasi-isometric.
  
  The image of a typical word with all $e_n$ and $f_n$ positive is shown in Figure \ref{fig:typicalm},
  corresponding to a series of rightward rotations at nodes distance one from the sides of the tree.

 {\bf Case $m > 0$:}
 
 In this case, we start with the same tree pair diagram for the $m=0$ case and apply
 $x_0$ on the right  $m$ times.  Each application of $x_0$ will change the negative tree by moving the
 root caret to a right caret and the topmost left caret to the root, if there is a left caret.
 If there is no left caret, a new caret will need to be added for each such application.
For each application of $x_0$ which requires a new caret, in the negative tree, that new caret will become the root caret and in the positive tree,
 the new caret will be added as the left child of the leftmost caret.
 Since there are $j_l$ left carets,  if $m \leq  j_l$, we do not need to add any
 carets and  the number of carets is
  $i_k+j_l+2+ \sum_{n=1}^k (| e_n|+1) + \sum_{n=1}^l (| f_n|+1) $ as before.  If $m> j_l$,
  we will need to add $m-j_l$ new carets and will have $i_k+j_l+2+ \sum_{n=1}^k (| e_n|+1) + \sum_{n=1}^l (| f_n|+1) + m-j_l $ carets.
  Again, these quantities give lengths which are comparable to word length in $\Z \wr \Z$.

 {\bf Case $m < 0$:}
 
 Again, we start with the same tree pair diagram for the $m=0$ case and apply
 $x_0^{-1}$ on the right $m$ times.  Each application of $x_0^{-1}$ will change the negative tree by moving the
 root caret to become a left caret and the topmost right caret to the root, if there is a right caret.
 If there is no right caret, a new caret will need to be added for each such application.
 Since there are $i_k$ left carets,  if $-m \leq  i_k$, we have that the number of carets is
  $i_k+j_l+2+ \sum_{n=1}^k (| e_n|+1) + \sum_{n=1}^l (| f_n|+1) $ as before.  If $-m> i_k$,
  we have  $i_k+j_l+2+ \sum_{n=1}^k (| e_n|+1) + \sum_{n=1}^l (| f_n|+1) -m + i_k $ carets.
  
  Thus in all cases $\phi$ does not distort distances more than linearly, so the
  subgroup $H$ isomorphic to $\Z \wr \Z$ is undistorted in $F$.
  
\qed.

We can obtain more precise estimates of the quasi-isometry constants using Fordham's
method \cite{blakegd} for computing exact lengths in $F$.  We can keep
track of the particular caret pairings and their weights and we find that the caret pairings
that occur are easily computed.  Caret pairing types are described in \cite{blakegd} and \cite{ct2}.
For example, in the case
where $m=0$ and both $l$ and $k$ are positive, we note that we have the following caret pairs:

\begin{itemize}
\item One caret pair of type $(L_0,L_0)$ from the leftmost carets, contributing no weight.
\item $j_l$ caret pairs of type $(L_L,L_L)$ from the left side and root, contributing weight $2 j_l$.
\item $i_k-1$ caret pairs of types $(R_*,R_*)$ not of type $(R_0,R_0)$, contributing weight $2 (i_k-1)$.
\item One caret pair of type $(R_0,R_0)$ from the rightmost carets, contributing no weight.
\item For each $e_n>0$, there will be a single pairing of type $(I_0,I_0)$ contributing weight 2
and $e_n-1$ pairings of type $(I_0,I_R)$, contributing weight $4 (e_n-1)$.
\item For each $e_n<0$, there will be a single pairing of type $(I_0,I_0)$ contributing weight 2
and $|e_n|-1$ pairings of type $(I_R,I_0)$, contributing weight $4 (|e_n|-1)$.
\item Similarly, for the interior carets from the left side of the tree, we have
for each $f_n$, there will be a single pairing of type $(I_0,I_0)$ contributing weight 2
and  $|f_n|-1$ pairings of type $(I_0,I_R)$ or $(I_R,I_0)$, contributing weight $4 (|f_n|-1)$.
\end{itemize}

These will give a total weight of $2 j_l+ 2 (i_k-1) + 2k + 4 \sum |e_n| + 2l + 4 \sum |f_n|$=
$2 j_l+ 2 i_k + 2k + 2l +  4 \sum |e_n| + 4 \sum |f_n| -2$ in the case when $m=0$, which compares
to the corresponding length in $\Z \wr \Z$  of $ 2 j_l +2 i_k+ \sum |e_n| +  \sum |f_n|$.

Again, these give lengths which are comparable to word length in $\Z \wr \Z$.  After similar
  analysis for other cases, we have  that for a word $w$ in $\Z \wr \Z$, we have:
  $$|w|_{\Z \wr \Z}-2 \leq |\phi(w)|_F \leq 4 |w|_{\Z \wr \Z}.$$

\section{$\Z \wr \Z$ as a subgroup of Baumslag's metabelian group}
\label{sec:metab}

Baumslag \cite{baumslag} introduced the group $G= <a,s,t | [s,t] , [a^t,a], a^s=a a^t>$ 
to show that a finitely presented metabelian group can contain free abelian subgroups
of infinite rank.   This group in fact contains $\Z \wr \Z$ -- all relators of the form
$[a^{t^i},a^{t^j}]$ are consequences of these three, so the subgroup $H$ generated
by $a$ and $t$ is isomorphic to $\Z \wr \Z$.

Here we examine the distortion of this subgroup in $G$

\begin{theorem}
The subgroup $H$ has at least exponential distortion in $G$.
\end{theorem}

{\bf Proof:}

First, we note that $s$ conjugates elements in $H$ to other elements in $H$ in
a manner illustrated here:

$$ a^{(s^2)}= (a^s)^s= (a a^t)^s = a^s (a^s)^t = a a^t (a a^t)^t = a a^t a^t a^{t^2} = a_0 a_1^2 a_2 $$

In terms of the notation described above, we have $a_n^s= a_n a_{n+1}$.
Further conjugation by $s$ leads to increasingly long words:
$$ a^{(s^3)} = (a_0 a_1^2 a_2)^s = a_0 a_1 a_1^2 a_2^2 a_2 a_3 = a_0 a_1^3 a_2^3 a_3$$
and we notice the occurrence of the binomial coefficients with repeated iteration.

In general, we have:

$$a^{s^n}= a_0^{n \choose 0} a_1^{n \choose 1}  \ldots a_n^{n \choose n}.$$

As an element of $G$, this has length $2n+1$ and it lies in the subgroup $H$, as there
is a representative with no occurrences of $s$.

To compute the length of this element in the subgroup $H$ with respect to
its generators $a$ and $t$, we use the method described in Section 2 and get

$$ |a^{s^n}|_H = 2n+ \sum_{i=0}^n {n \choose i}= 2n+2^{n}.$$

Thus we have $ |a^{s^n}|_H=2n+2^{n}$ while $|a^{s^n}|_G=2n+1$, 
so the wreath product $\Z \wr \Z$ is exponentially distorted in $G$.

\qed

\bibliographystyle{plain}

\def\cprime{$'$}


%

\end{document}